\documentclass[12pt]{article}
\usepackage{graphicx}
\usepackage{color}
\catcode`\@=11 \@addtoreset{equation}{section}

\catcode`\@=12

\newtheorem{Theorem}{Theorem}[section]
\newtheorem{Proposition}{Proposition}[section]
\newtheorem{Lemma}{Lemma}[section]
\newtheorem{Corollary}{Corollary}[section]

\newcommand{\bTheorem}[1]{
\begin{Theorem} \label{T#1} }
\newcommand{\eT}{\end{Theorem}}

\newcommand{\bProposition}[1]{
\begin{Proposition} \label{P#1}}
\newcommand{\eP}{\end{Proposition}}

\newcommand{\bLemma}[1]{
\begin{Lemma} \label{L#1} }
\newcommand{\eL}{\end{Lemma}}

\newcommand{\bCorollary}[1]{
\begin{Corollary} \label{C#1} }
\newcommand{\eC}{\end{Corollary}}

\newcommand{\bFormula}[1]{
\begin{equation} \label{#1}}
\newcommand{\eF}{\end{equation}}

\newcommand{\Ov}[1]{\overline{#1}}

\newcommand{\DC}{C^\infty_c}

\newcommand{\vr}{\varrho}
\newcommand{\vre}{\vr_\ep}

\newcommand{\vue}{\vu_\ep}

\newcommand{\vu}{\vc{u}}
\newcommand{\vc}[1]{{\bf #1}}

\newcommand{\Div}{{\rm div}_x}
\newcommand{\Grad}{\nabla_x}

\newcommand{\tn}[1]{\mbox {\F #1}}
\newcommand{\dx}{{\rm d} {x}}
\newcommand{\dt}{{\rm d} t }

\newcommand{\intO}[1]{\int_{\Omega} #1 \ \dx}

\newcommand{\ep}{\varepsilon}

\font\F=msbm10 scaled 1000

\definecolor{Cgrey}{rgb}{0.85,0.85,0.85}
\definecolor{Cblue}{rgb}{0.50,0.85,0.85}
\definecolor{Cred}{rgb}{1,0,0}
\definecolor{fancy}{rgb}{0.10,0.85,0.10}

\newcommand\Cbox[2]{%
    \newbox\contentbox%
    \newbox\bkgdbox%
    \setbox\contentbox\hbox to \hsize{%
        \vtop{
            \kern\columnsep
            \hbox to \hsize{%
                \kern\columnsep%
                \advance\hsize by -2\columnsep%
                \setlength{\textwidth}{\hsize}%
                \vbox{
                    \parskip=\baselineskip
                    \parindent=0bp
                    #2
                }%
                \kern\columnsep%
            }%
            \kern\columnsep%
        }%
    }%
    \setbox\bkgdbox\vbox{
        \color{#1}
        \hrule width  \wd\contentbox %
               height \ht\contentbox %
               depth  \dp\contentbox
        \color{black}
    }%
    \wd\bkgdbox=0bp%
    \vbox{\hbox to \hsize{\box\bkgdbox\box\contentbox}}%
    \vskip\baselineskip%
}

\date{}

\begin{document}


\title{Multiple scales and singular limits for compressible rotating fluids with general initial data}

\author{Eduard Feireisl\thanks{Eduard Feireisl acknowledges the support of the project LL1202 in the
programme ERC-CZ funded
by the Ministry of Education, Youth and Sports of the Czech Republic.} \and
Anton\' \i n Novotn\' y }

\maketitle

\bigskip

\centerline{Institute of Mathematics of the Academy of Sciences of the Czech Republic}

\centerline{\v Zitn\' a 25, 115 67 Praha 1, Czech Republic}

\centerline{Charles University in Prague, Faculty of Mathematics and Physics, Mathematical Institute}

\centerline{Sokolovsk\' a 83, 186 75 Praha 8,
Czech Republic}

\centerline{and}

\centerline{IMATH, EA 2134, Universit\' e du Sud Toulon-Var
BP 20132, 83957 La Garde, France}


\maketitle

\bigskip





\begin{abstract}

We study the singular limit of a rotating compressible fluid described by a scaled barotropic Navier-Stokes system, where the
Rossby number $= \ep$, the Mach number $ =
\ep^m$, the {Reynolds number} $= \ep^{-\alpha}$, and the \emph{Froude number} $= \ep^n$ are proportional to a small parameter $\ep \to 0$.
The inviscid planar Euler system is identified as the limit problem. The proof is based on the application of the method of relative entropies
and careful analysis of oscillatory integrals describing the propagation of Rossby-acoustic waves.

\end{abstract}


\section{Introduction}
\label{i}

The basic system of equations considered in this paper are the \emph{continuity equation} and the \emph{momentum equation} describing the time evolution of the
mass density $\vr = \vr(t,x)$ and the (relative) velocity $\vu = \vu(t,x)$ of a compressible, rotating fluid:

\Cbox{Cgrey}{

\bFormula{i1}
\partial_t \vr + \Div (\vr \vu) = 0,
\eF
\bFormula{i2}
\partial_t (\vr \vu) + \Div (\vr \vu \otimes \vu) + \frac{1}{\ep} \vr (\vc{f} \times \vu) + \frac{1}{\ep^{2m}} \Grad p(\vr) =
\ep^\alpha \Div \tn{S}(\Grad \vu) + \frac{1}{\ep^{2n}} \vr \Grad G,
\eF
\bFormula{i3}
\tn{S}(\Grad \vu) = \mu \left( \Grad \vu + \Grad^t \vu - \frac{2}{3} \Div \vu \tn{I} \right) + \eta \Div \vu \tn{I}, \ \mu > 0, \ \eta \geq 0.
\eF

}

\noindent
The fluid is confined to an infinite slab
\bFormula{i4}
\Omega = R^2 \times (0,1),
\eF
where it satisfies the \emph{slip condition}
\bFormula{i5}
\vu \cdot \vc{n}|_{\partial \Omega} = [\tn{S}(\Grad \vu) \cdot \vc{n}]_{\rm tan}|_{\partial \Omega} = 0
\eF
imposed on the horizontal boundary.

The model (\ref{i1} - \ref{i5}) may be viewed as a crude approximation ($f$-plane model) of the Earth atmosphere in a plane tangent to the Earth at a certain
latitude, see Vallis \cite[Chapter 2, Section 2.3]{VALL}. Accordingly, the gravitational force is taken parallel to the vertical projection of the rotation axis:
 \[
\vc{f} = [0,0,1], \ \Grad {G} = [0,0,-1],
\]

The momentum equation (\ref{i2}) contains a small parameter $\ep$
related to different \emph{characteristic numbers} resulting from
the scale analysis: \emph{Rossby number} $= \ep$, \emph{Mach
number} $ = \ep^m$, \emph{Reynolds number} $= \ep^{-\alpha}$,
\emph{Froude number} $= \ep^n$, see Klein \cite{Klein}.
We consider the singular limit
problem for $\ep \searrow 0$ in the multiscale regime:
\bFormula{i8} \frac{m}{2} > n \geq 1, \ \alpha
> 0 \eF
for the
\emph{ill-prepared initial data}
\bFormula{i6} \vr(0, \cdot) =
\vr_{0,\ep} = \tilde \vr_\ep + \ep^m \vr^{(1)}_\ep ,\ \vu(0,
\cdot) = \vu_0, \eF
where $\tilde \vre$ is a solution to the
static problem
\bFormula{i7} \Grad p(\tilde \vre) = \ep^{2(m-n)} \tilde \vre
\Grad G \ \mbox{in}\ \Omega.\eF
In particular, since $n \geq 1$, the centrifugal force, parallel
to the vertical axis, is dominated by gravitation (cf. Durran
\cite{Dur1}).

Formally, it is not difficult to identify the limit problem. Indeed fast rotation is expected to eliminate the vertical motion, the vanishing  viscosity
(high Reynolds number) makes the limit system \emph{inviscid} (hyperbolic), while the low Mach number regime drives the fluid to \emph{incompressibility}. The limit problem is
therefore expected to be the incompressible Euler system for the planar velocity field $\vc{v} = [v_1,v_2]$,
\bFormula{i9}
\partial_t \vc{v} + \Div (\vc{v}  \otimes \vc{v} ) + \Grad \Pi = 0, \ \Div \vc{v} = 0 \ \mbox{in} \ (0,T) \times R^2.
\eF
Our main goal is to put these formal arguments on rigorous grounds.

The phenomena discussed above have been investigated by many authors. The fact that highly rotating fluids become planar (two-dimensional), and, accordingly,
fast rotation has a regularizing effect, was observed by Babin, Mahalov, and Nicolaenko \cite{BaMaNi2}, \cite{BaMaNi1}, Bresch, Desjardins, and Gerard-Varet, \cite{BrDeGV}, Chemin et al. \cite{CDGG}, among others. The inviscid limit is a well studied and partially still open challenging problem, see
Clopeau, Mikeli´c, Robert \cite{CloMikRob}, Kato \cite{Kato}, Masmoudi  \cite{MAS6}, \cite{MAS7}, \cite{MAS1}, Sammartino and Caflisch \cite{SamCaf1},
\cite{SamCaf2}, Swann \cite{Swa}, Temam and Wang \cite{TemWan1}, \cite{TemWan2}, to name only a few. The low Mach number limits were proposed in the pioneering papers by Ebin \cite{EB1}, and Klainerman and Majda \cite{KM1}, and later reexamined in the context of weak solutions by Lions and Masmoudi \cite{LIMA1}, see also the
survey by Danchin \cite{DANC4}, Gallagher \cite{Gallag}, and Schochet \cite{SCH2}. To the best of our knowledge, the simultaneous effects of these three
mechanisms has not yet been treated in the literature.

The present paper may be viewed as complementary to our previous study \cite{FeNo7}, where we examined the ``single-scale'' limit corresponding to
\bFormula{i10}
n = 0, \ m = 1, \ \alpha > 0.
\eF
Although this problem looks formally very similar to the present setting, the methods employed as well as the limit system are \emph{different}, cf.
\cite{FeNo7}. The central issue
to be discussed is the behavior of the oscillatory part of solutions to the scaled system. These are described in the case (\ref{i10}) by a system of linear equations with \emph{constant} coefficients, while the situation (\ref{i8}) gives rise to a problem with \emph{coefficients depending on the scaling parameter} $\ep$.

Similarly to \cite{FeNo7}, our approach is based on the concept of finite energy weak solutions satisfying the \emph{relative entropy inequality} identified in
\cite{FeNoSun}, see Section \ref{p}. After collecting the necessary preliminary material, we state the main result in Section \ref{m}. Section \ref{e} reviews the basic estimates, independent of $\ep \searrow 0$, for solutions of the family of scaled problems. Section \ref{d} is the heart of the paper. We establish  decay estimates on the oscillatory part of solutions to the scaled problems by means of careful analysis of certain oscillatory integrals. Here, inspired by the analysis of Guo, Peng, and Wang \cite{GuPeWa}, we make use of frequency cut-off operators and estimates based on van Corput's lemma. The proof of convergence towards the limit system is completed in Section \ref{c}.

\section{Preliminaries, weak solutions, relative entropy inequality}
\label{p}

We suppose that the \emph{pressure} $p \in C[0, \infty) \cap
C^3(0, \infty)$ is a given function of the density enjoying the
following property \bFormula{p1} p(0) = 0, \ p'(\vr) > 0 \
\mbox{for all}\ \vr > 0, \ \lim_{\vr \to \infty}
\frac{p'(\vr)}{\vr^{\gamma - 1}} , \ \gamma > \frac{3}{2}. \eF In
addition, without loss of
generality, we assume that $p$ is ``normalized'' setting
\bFormula{p1a} p'(1) = 1. \eF

Finally, we introduce the pressure potential $H$,
\bFormula{p1b}
H(\vr) = \vr \int_1^\vr \frac{p(z)} {z^2} \ {\rm d}z,
\eF
noting that
\[
H''(\vr) = \frac{p'(\vr)}{\vr}, \ H''(1) = 1.
\]

\subsection{Static solutions}

As already mentioned above, the static solutions $\tilde \vre$
solve the problem (\ref{i7}), specifically, we take
\bFormula{p2} H'(\tilde \vre) = \ep^{2(m-n)}G +
H'(1),\;\mbox{where $G(x)=-x_3$}; \eF
whence  \bFormula{p3} \tilde\vre(x)=\tilde\vre(x_3),\;\sup_{x_3
\in [0,1]} | \tilde \vre(x_3) - 1 | \leq c \ep^{2(m-n)}. \eF

As indicated by our choice of the initial data (\ref{i6}), the solutions of the evolutionary problem (\ref{i1} - \ref{i3}), (\ref{i5}), (\ref{i6}) satisfy \emph{far field conditions}
in the form
\bFormula{p4}
\vr \to \tilde \vre , \ \vu \to 0 \ \mbox{as}\ |x| \to \infty.
\eF

\subsection{Finite energy weak solutions}

\label{few}

We say that $[\vr, \vu]$ is a \emph{finite energy weak solution} of the problem (\ref{i1} - \ref{i3}), (\ref{i5}), (\ref{i6}), (\ref{p4})
on the space-time cylinder $(0,T) \times \Omega$ if the following holds:
\begin{itemize}
\item {\bf Regularity.}
The functions $\vr$, $\vu$ belong to the class
\bFormula{few1}
\vr \geq 0,\ (\vr - \tilde \vre) \in L^\infty(0,T; L^2 + L^\gamma (\Omega)), \ \vu \in L^2(0,T; W^{1,2}(\Omega;R^3)), \ \vu \cdot \vc{n} = u_3|_{\partial \Omega} = 0.
\eF
\item {\bf Equations.} The equation of continuity (\ref{i1}) and the momentum equation (\ref{i2}) are replaced by integral identities
\bFormula{few2} \int_0^T \intO{ \left( \vr \partial_t \varphi +
\vr \vu \cdot \Grad \varphi \right) } \ \dt = - \intO{ \vr_{0,
\ep} \varphi(0, \cdot) } \eF for any $\varphi \in \DC([0,T) \times
\Ov{\Omega})$, and \bFormula{few3} \int_0^T \intO{ \left( \vre
\vue \cdot \partial_t \varphi + (\vr \vu \otimes \vu): \Grad
\varphi - \frac{1}{\ep} \vr (\vc{f} \times \vu){\cdot \varphi}
+ \frac{1}{\ep^{2m}} p(\vr) \Div \varphi \right) } \ \dt \eF
\[
= \int_0^T \intO{ \left( \ep^\alpha \tn{S}(\Grad \vu) : \Grad
\varphi - \frac{1}{\ep^{2n}} \vr \Grad G \cdot \varphi \right) } \
\dt - \intO{ {\vr}_{0, \ep} \vu_{0, \ep} \cdot \varphi (0,
\cdot) }
\]
for any $\varphi \in \DC([0,T) \times \Ov{\Omega}; R^3)$, $\varphi \cdot \vc{n}|_{\partial \Omega}  = 0$.

\item {\bf Energy.}
The energy inequality
\bFormula{few4}
\intO{ \left[ \frac{1}{2} \vr |\vu|^2 + \frac{1}{\ep^{2m}} \left( H(\vr) - H'(\tilde \vre)(\vr - \tilde \vre) - H(\tilde \vre) \right) \right](\tau, \cdot) }
+ \ep^\alpha \int_0^\tau \intO{ \tn{S}(\Grad \vu) : \Grad \vu } \ \dt
\eF
\[
\leq \intO{ \left[ \frac{1}{2} \vr_{0,\ep} |\vu_{0,\ep}|^2 + \frac{1}{\ep^{2m}} \left( H(\vr_{0,\ep}) - H'(\tilde \vre)(\vr_{0,\ep} - \tilde \vre) - H(\tilde \vre) \right) \right] }
\]
holds for a.a. $\tau \in (0,T)$.

\end{itemize}

Note that the \emph{existence theory} in the class of finite energy weak solutions was developed by Lions \cite{LI4} and later extended in \cite{FNP}
to the sofar ``critical'' adiabatic exponent $\gamma > \frac{3}{2}$.

\subsection{Relative entropy}

For future analysis, it is convenient to replace the energy
inequality (\ref{few4}) by the \emph{relative entropy inequality}
containing more transparent piece of information on the
asymptotic behavior of solutions for $\ep \to 0$. To this end, we
introduce the \emph{relative entropy} functional \bFormula{p5}
\mathcal{E}_\ep \left( \vr, \vu \Big| r, \vc{U} \right) = \intO{
\left[ \frac{1}{2} \vr |\vu - \vc{U}|^2 + \frac{1}{\ep^{2m}} \Big(
H(\vr) - H'(r)(\vr - r) - H(r) \Big) \right]}, \eF cf.
\cite{FeJiNo}, \cite{FeNoSun}, Germain \cite{Ger}.

It can be shown, see \cite{FeJiNo}, that \emph{any} finite energy weak solution $[\vr, \vu]$
specified in Section \ref{few} satisfies the
\emph{relative entropy inequality}:
\bFormula{p6}
\mathcal{E}_\ep \left( \vr, \vu \ \Big| \ r, \vc{U} \right) (\tau)
+ \ep^\alpha \int_0^\tau \intO{ \Big( \tn{S} (\Grad \vu) - \tn{S} (\Grad \vc{U}) \Big) : \Big( \Grad \vu - \Grad \vc{U} \Big) } \ \dt \leq
\eF
\[
\mathcal{E}_\ep \left( \vr_{0,\ep}, \vu_{0,\ep} \ \Big| \ r(0,\cdot) , \vc{U}(0,\cdot) \right)
\]
\[
+ \int_0^\tau \intO{  \vr \left( \partial_t \vc{U} + \vu \cdot \Grad \vc{U} \right) \cdot \left( \vc{U} - \vu \right) } \ \dt
\]
\[
+ \ep^\alpha \int_0^\tau \intO{ \tn{S} (\Grad \vc{U}) : \Grad (\vc{U} - \vu ) } \ \dt + \frac{1}{\ep} \int_0^\tau \intO{
\vr (\vc{f} \times \vu ) \cdot (\vc{U} - \vu) } \ \dt
\]
\[
+ \frac{1}{\ep^{2m}} \int_0^\tau \intO{ \Big[ (r - \vr) \partial_t
H'(r) + { \Grad\Big(  H'(r)-H'(\tilde\vr_\ep)\Big) } \cdot (r
\vc{U} - \vr \vu ) \Big] } \ \dt
\]
\[
- \frac{1}{\ep^{2m}} \int_0^\tau \intO{ \Div \vc{U} \Big( p(\vr) -
p(r) \Big) } \ \dt { - \frac{1}{\ep^{2n}} \int_0^\tau \intO{
(\vr-r) \Grad G \cdot \vc{U} } \ \dt}
\]
for all sufficiently smooth ``test functions'' $r$, $\vc{U}$
satisfying
\[
\vc{U} \cdot \vc{n}|_{\partial \Omega} = 0,\ r > 0,\  \vc{U}, \ (r - \tilde \vre) \ \mbox{compactly supported in}\ \Ov{\Omega}.
\]
Note the assumptions concerning the decay and regularity  can be relaxed to basically any couple $[r, \vc{U}]$ for which (\ref{p6}) makes sense
via the standard density argument, see \cite{FeJiNo}.

\section{Main result}
\label{m}

For a vector field $\vc{b}=[b_1,b_2,b_2]$, we introduce the horizontal component
$\vc{b}_h= [b_1,b_2]$ writing $\vc{b} = [\vc{b}_h, b_3]$. Similarly, we use the symbols
$\nabla_h$, ${\rm div}_h$ to denote the differential operators acting on the horizontal variables only. Finally, the symbol
$\vc{H}$ denotes the Helmholtz projection onto the space of solenoidal functions in $\Omega$, while $\vc{H}_h$ is the Helmholtz projection in $R^2$.

{Let
$$
\vc v_0\in W^{k,2}(R^2;R^2),\;k\ge 3,\;{\rm div}\vc v_0=0,
$$
be given.
It is well known (see for instance Kato and Lai \cite{KaLai}) that the Euler system
(\ref{i9}) supplemented with the initial data
$$
\vc v(0)=\vc v_0
$$
admits a regular solution $(\vc v,\Pi)$, unique in the class
$$
\vc v\in C([0,T]; W^{k,2}(R^2;R^2)),\;\partial_t\vc v\in  C([0,T];
W^{k-1,2}(R^2;R^2)),\;
 \Pi\in C([0,T]; W^{k,2}(R^2)).
$$
}

 We are ready to formulate our main result.

\Cbox{Cgrey}{

\bTheorem{m1}

Let the pressure $p=p(\vr)$ satisfy the hypotheses (\ref{p1}), (\ref{p1a}).
Suppose that the exponents $\alpha$, $m$, $n$ are given such that
\[
\alpha > 0, \ \frac{m}{2} > n \geq 1.
\]
Let the initial data $[\vr_{0,\ep}, \vu_{0, \ep}]$ be given by (\ref{i6}), where
the stationary states $\tilde \vre$ satisfy (\ref{p2}),
\bFormula{m1}
\| \vr^{(1)}_{0, \ep} \|_{L^2 \cap L^\infty (\Omega)} \leq c,\ \vr^{(1)}_{0,\ep} \to \vr^{(1)}_0 \ \mbox{in}\ L^2(\Omega),\
\vu_{0,\ep} \to \vu_0 \ \mbox{in}\ L^2(\Omega;R^3),
\eF
with
\bFormula{m2}
\vr^{(1)}_0 \in W^{k-1,2}(\Omega),\
\vu_0 \in W^{k,2}(\Omega;R^3) \ \mbox{for a certain}\ k \geq 3.
\eF
Let $[\vre,\vue]$ be a finite energy weak solution of the problem (\ref{i1} - \ref{i3}), (\ref{i5}), (\ref{i6}), (\ref{p4}) in the
space-time cylinder
$(0,T) \times \Omega$.

Then
\bFormula{m3}
{\rm ess} \sup_{t \in (0,T)} \| \vre (t, \cdot) - \tilde \vre \|_{(L^2 + L^\gamma)(\Omega)} \leq \ep^m c
\eF
\bFormula{m4}
\sqrt{\vre} \vue \to \vc{v} \left\{ \begin{array}{l} \mbox{weakly-(*) in} \ L^\infty(0,T;
L^2  (\Omega; R^3)), \\ \\ \mbox{strongly in}\ L^1_{\rm loc}((0,T) \times {\Omega};R^3),
\end{array} \right.
\eF
where $\vc{v} = [\vc{v}_h,0]$ is the unique solution of the Euler system (\ref{i9}), with the initial data
\[
\vc{v}(0, \cdot) = \vc{H}_h \left[ \int_0^1 \vu_0(x_h, x_3) \ {\rm d}x_3 \right].
\]
\eT

}

The rest of the paper is devoted to the proof of Theorem \ref{Tm1}.

\section{Uniform bounds}
\label{e}

We start with the nowadays standard estimates that follow directly from the energy inequality (\ref{few4}).
These are conveniently formulated in terms of
the \emph{essential} and \emph{residual} components of a measurable function $h$,
\[
h = h_{\rm ess} + h_{\rm res},
\]
\[
h_{\rm ess} = \chi(\vre) h,\ \chi \in \DC(0, \infty), \ 0 \leq \chi \leq 1,\
\chi = 1 \ \mbox{in an open interval contaning} \  1,
\]
\[
h_{\rm res} = (1 - \chi(\vre)) h.
\]

Since the initial data are given by (\ref{i6}), where the
functions $\vr^{(1)}_{0, \ep}$, $\vu_{0,\ep}$ satisfy the
hypotheses (\ref{m1}), (\ref{m2}), the integral on the {
right-hand} side of (\ref{few}) is bounded uniformly for $\ep
\searrow 0$. As the stationary states $\tilde \vre$ are chosen to
satisfy (\ref{p3}), we deduce the following bounds independent of
$\ep$: \bFormula{rr2} {\rm ess} \sup_{t \in (0,T)} \left\|
\sqrt{\vre} \vue \right\|_{L^2(\Omega;R^3)} \leq c, \eF
\bFormula{rr3} {\rm ess} \sup_{t \in (0,T)} \left\| \left[
\frac{\vre - \tilde \vre}{\ep^m} \right]_{\rm ess}
\right\|_{L^2(\Omega)} \leq c, \eF \bFormula{rr4} {\rm ess}
\sup_{t \in (0,T)} \left\| [\vre ]_{\rm res}
\right\|^\gamma_{L^\gamma(\Omega)} + {\rm ess} \sup_{t \in (0,T)}
\left\| [ 1 ]_{\rm res}  \right\|_{L^1(\Omega)} \leq \ep^{2m} c,
\eF and \bFormula{rr5} \ep^\alpha \int_0^T \intO{ \left| \Grad
\vue + \Grad \vue^t - \frac{2}{3} \Div \vue \tn{I} \right|^2 } \
\dt \leq c, \eF cf. \cite[Section 2]{FeGaNo}.

Obviously, the estimates (\ref{rr2}), (\ref{rr3}) yield (\ref{m3}), which, combined with (\ref{p3}) gives rise to
\bFormula{e1}
\vre \to 1 \ \mbox{in}\ L^\infty(0,T; L^q_{\rm loc}(\Omega)) \ \mbox{for any} \ 1 \leq q < \gamma.
\eF
Moreover, combining (\ref{e1}) with (\ref{rr2} - \ref{rr4}) we obtain
\bFormula{e2}
\sqrt{\vre} \vue \to \vu \ \mbox{weakly-(*) in} \ L^\infty(0,T; L^2(\Omega;R^3)),
\eF
and
\bFormula{e3}
\vre \vue \to \vu \ \mbox{weakly-(*) in}\ L^\infty(0,T; L^2 + L^{2 \gamma /(\gamma + 1)}(\Omega;R^3)),
\eF
passing to suitable subsequences as the case may be.

Finally, we may let $\ep \to 0$ in the equation of continuity to deduce that
\bFormula{e4}
\Div \vu = 0, \ \vu \cdot \vc{n}|_{\partial \Omega} = 0 \ \mbox{in the sense of distributions in} \ (0,T) \times \Omega.
\eF

\section{Decay estimates and oscillatory integrals}
\label{d}

With our convention (\ref{p1a}),
the equation describing the oscillatory part of solutions  reads
\bFormula{d1}
\ep^m \partial_t s + \Div \vc{V} = 0,
\eF
\bFormula{d2}
\ep^m \partial_t \vc{V} + \omega \vc{f} \times \vc{V} + \Grad s = 0, \ \omega = \ep^{m-1}, \ \vc{V} \cdot \vc{n}|_{\partial \Omega} = 0,
\eF
cf. \cite{FeNo7}.
Re-scaling in the time we arrive at
\bFormula{d3}
\partial_t s + \Div \vc{V} = 0,
\eF
\bFormula{d4}
\partial_t \vc{V} + \omega \vc{b} \times \vc{V} + \Grad s = 0, \ \vc{V} \cdot \vc{n}|_{\partial \Omega} = 0,
\eF
with the operator
\[
\mathcal{B}(\omega): \left[ \begin{array}{c} s \\ \vc{V} \end{array} \right] \mapsto \left[ \begin{array}{c} \Div \vc{V} \\ \omega \vc{f} \times \vc{V} + \Grad s \end{array} \right].
\]

The operator $\mathcal{B}$ is skew symmetric in the space $L^2(\Omega) \times L^2(\Omega;R^3)$, with the domain of definition
\[
\mathcal{D}[\mathcal{B}(\omega)] = \left\{ [r, \vc{V}] \ \Big| \ r \in W^{1,2}(\Omega), \vc{V} \in L^2(\Omega;R^3), \Div \vc{V} \in L^2(\Omega),
\vc{V} \cdot \vc{n} = V_3|_{\partial \Omega} = 0
\right\},
\]
and the kernel
\[
\mathcal{N}(\mathcal{B}(\omega)) = \left\{ [q, \vc{v}] \ \Big| \ q = q(x_h), \ q \in W^{1,2}(R^2),
\ \vc{v} = [ \vc{v}_h (x_h), 0],\ {\rm div}_h \vc{v}_h = 0 , \
{\omega} \vc{f} \times \vc{v} + \Grad q = 0 \right\}.
\]

Let $\mathcal{P}(\omega)$ denote the projection
\[
\mathcal{P}(\omega) : L^2(\Omega) \times L^2(\Omega;R^3) \to \mathcal{N}(\mathcal{B}(\omega)).
\]
Exactly as in \cite[Section 4.1.1]{FeNo7} we can show that
\[
\mathcal{P}(\omega)[r, \vc{U}] = [q , \vc{v}]
\]
if
\bFormula{pro}
- \Delta_h q + \omega^2 q = \omega \int_0^1 {\rm curl}_h \vc{U}_h \ {\rm d}x_3 + \omega^2 \int_0^1 r \ {\rm d}x_3 \ \mbox{in}\ R^2, \ \vc{v}=[v_1,v_2,0],\
\omega v_1 = - \partial_{x_2}q, \ \omega v_2 = \partial_{x_1} q.
\eF

\subsection{Spectral analysis}

Thanks to our special choice of the geometry of the spatial domain $\Omega$, we may
reformulate
the problem (\ref{d3}), (\ref{d4}) in terms of the Fourier variables. To this end, we observe, exactly as in \cite{FeNo7}, that the underlying spatial domain
$\Omega$ may be equivalently replaced by
\[
\tilde \Omega = R^2 \times \mathcal{T}^1 ,
\]
where
\[
\mathcal{T}^1 = [-1,1]_{\{ -1, 1 \} }
\]
is a ``flat'' sphere, and where
$
s, \ \vc{V}_h \ \mbox{were extended as even functions of the vertical variable}\ x_3,
$
while
$
V_3 \ \mbox{was extended as odd in} \ x_3.
$
In other words, all quantities are understood as $2-$periodic functions in the vertical $x_3$ variable.

Accordingly, for each function $g \in L^2(\tilde \Omega)$, we introduce its Fourier representation
\[
\hat g (\xi, k), \ \xi = [\xi_1, \xi_2] \in R^2, \ k \in Z,
\]
where
\[
\hat g(\xi,k) = \frac{1}{\sqrt{2}} \int_{-1}^1 \int_{R^2} \exp
\left(-{\rm i} \xi \cdot x_h \right) g(x_h,x_3)  \exp \left(- {\rm
i} k  x_3 \right){\ {\rm d}x_h}\ {\rm d}x_3.
\]
We have
\[
g(x_h, x_3) = \sum_{k \in Z} \mathcal{F}^{-1}_{\xi \to x_h} \left[ \hat g(\xi,k ) \right] \exp \left({\rm i}k x_3 \right),
\]
where the symbol $\mathcal{F}_{x_h \to \xi}$ denotes the standard Fourier transform on $R^2$.

Thus the problem (\ref{d3}), (\ref{d4}) can be written in the form
\bFormula{d5}
\frac{{\rm d}}{{\rm d}t} \left[ \begin{array}{c} \hat s(t, \xi, k) \\ \hat {V}_1(t, \xi, k)
\\ \hat {V}_2(t, \xi, k) \\ \hat {V}_3(t, \xi, k)
 \end{array}
\right] + {\rm i} \left[ \begin{array}{cccc}
                   0 & \xi_1 & \xi_2 & k \\
                   \xi_1& 0 & \omega{\rm i} & 0 \\
                   \xi_2& - \omega{\rm i} & 0 & 0 \\
                    k& 0 & 0 & 0
                 \end{array}
 \right] \left[ \begin{array}{c} \hat s(t, \xi, k) \\ \hat {V}_1(t, \xi, k)
\\ \hat {V}_2(t, \xi, k) \\ \hat {V}_3(t, \xi, k)
 \end{array}
\right] = 0 , \left[ \begin{array}{c} \hat s(0, \xi, k) \\ \hat \vc{V}(0, \xi, k) \end{array}
\right] = \left[ \begin{array}{c} \hat s_0(\xi, k) \\ \hat \vc{V}_0(\xi, k) \end{array}
\right];
\eF
with the Hermitian matrix
\[
\mathcal{A}(\xi,k, \omega) = \left[ \begin{array}{cccc}
                   0 & \xi_1 & \xi_2 & k \\
                   \xi_1& 0 & \omega {\rm i} & 0 \\
                   \xi_2& - \omega {\rm i} & 0 & 0 \\
                    k& 0 & 0 & 0
                 \end{array}
 \right].
\]
Of course, solutions of (\ref{d5}) depend also on the parameter $\omega = \ep^{m-1}$.

\subsubsection{Spectral properties of the matrix $\mathcal{A}$}

After a bit tedious but straightforward manipulation (see \cite{FeGaNo}), we can check that $\mathcal{A}(\xi,k,\omega)$
possesses four eigenvalues
\bFormula{d6}
\begin{array}{c}
\lambda_{1} ( {|\xi|^2 ,k, \omega) = \left[ \frac{ \omega^2 + | \xi |^2
+ k^2 + \sqrt{ (\omega^2 + | \xi |^2 + k^2)^2 - 4 \omega^2 k^2
}}{2} \right]^{1/2} ,\ \lambda_2 (|\xi|^2, k, \omega)  = - \lambda_1
(|\xi|^2}, k, \omega),
\\ \\
\lambda_{3} (|\xi|^2,k, \omega) = \left[ \frac{ \omega^2 + | \xi |^2 +
k^2 - \sqrt{ (\omega^2 + | \xi |^2 + k^2)^2 - 4 \omega^2 k^2 }}{2}
\right]^{1/2},\ \lambda_4(|\xi|^2, k, \omega) = - \lambda_3 (|\xi|^2, k,
\omega).
\end{array}
\eF Note that $\lambda_3(|\xi|^2, 0 ,\omega) = \lambda_4(|\xi|^2,0,
\omega) = 0$ are the zero eigenvalues corresponding to the
non-void kernel of the matrix $\mathcal{A}$- the {Fourier
image of the} null-space of the operator $\mathcal{B}(\omega)$,
see \cite{FeGaNo}.

As for the eigenvectors $[q,v_1,v_2,v_3]$, we have \bFormula{ei3}
\xi_1 v_1 + \xi_2 v_2 +k v_3 = \lambda q, \ \xi_1 q + {\rm
i}\omega v_2 = \lambda v_1, \ \xi_2 q - {\rm i} \omega v_1 =
\lambda v_2, k q = \lambda v_3, \eF from which we immediately
deduce \bFormula{ei1} v_1 = \mu (\lambda \xi_1 + {\rm i}\omega
\xi_2), \ v_2 = \mu (\lambda \xi_2 - {\rm i}\omega \xi_1 ), \ \mu
\lambda |\xi|^2 = \lambda q - kv_3,\ ,\ k q = \lambda v_3, \eF
where $\mu$ is a free parameter that is fixed to normalize the
norm of the eigenvector to be one.

\subsubsection{Eigenvectors}
\label{EV}

{We denote by ${\bf E}={\bf E}(\xi,k,\omega) = [q,v_1,v_2,v_3]$ the normalized
eigenvectors. Our goal is to show that diagonalizing matrices $\mathcal{Q}$, $\mathcal{Q}^T$, formed by the eigenvectors, are $L^p-$multipliers
in the $\xi$ variable restricted to compact subsets of $R^2 \setminus \{0\}$. This amounts to showing that
\begin{equation}\label{multiplier}
\sup_{\omega\in (0,1)}\max_{0<a\le|\xi|\le
b<\infty}\Big|\nabla_\xi^A E_j(\xi,k,\omega)\Big|\le
c=c(A,a,b,k),\; j=1,2,3,4,\; k\in Z
\end{equation}
with any multi-index $A=(A_1,A_2)$.}

We distinguish
{two} cases:

\medskip

{\bf Case} { $\lambda= 0$:}

\medskip

In this case, we necessarily have $k=0$ and {
$\lambda=\pm\lambda_3$}, and the orthonormal basis of eigenvectors
can be taken in the form
\[
\vc{E}_1 = \mu \Big[ -{\rm i} \omega, - \xi_2 , \xi_1 , 0 \Big],\ \mu = \left( |\xi|^2 + \omega^2 \right)^{-1/2},\
\vc{E}_2 = [0,0,0,1].
\]

Clearly (\ref{multiplier}) holds.

\medskip

{\bf Case {$\lambda\neq 0$} :}

\medskip

We find that
\[
v_1 = \mu (\lambda \xi_1 + {\rm i}\omega \xi_2), \ v_2 = \mu (\lambda \xi_2 - {\rm i}\omega \xi_1 ),
\frac{\lambda^2 - k^2}{\lambda^2} q = \mu (\xi_1^2 +\xi_2^2 ), \ v_3 = \frac{k}{\lambda} q.
\]
Thus, the corresponding normalized eigenvector has the form
\[
\vc{E} = \mu \left[ \frac{\lambda^2 |\xi|^2}{\lambda^2 - k^2}, \lambda \xi_1
+ {\rm i}\omega \xi_2, \lambda \xi_2 - {\rm i}\omega \xi_1 ,
\frac{k \lambda |\xi|^2}{\lambda^2 - k^2}  \right],
\]
with
\[
\mu = \left[ \frac{\lambda^4 + |k|^2 \lambda^2}{(\lambda^2 - k^2)^2} |\xi|^4
+ (\lambda^2 + \omega^2) |\xi|^2  \right]^{-1/2}.
\]

We consider first the case $\lambda=\pm\lambda_1$. We  check that
{
\begin{equation}\label{ev} \lambda_1^2\ge|\xi|^2/2,\quad
(\omega^2+|\xi|^2+ k^2)^2-4\omega^2k^2\ge|\xi|^4,
\end{equation}
}
\[
\lambda_1^2 - k^2 = \frac{ \omega^2 - k^2 + |\xi|^2 +
\sqrt{(\omega^2 - k^2)^2  + |\xi|^2 \left(|\xi|^2 + 2 (\omega^2 +
k^2) \right) }} {2} \geq \frac{|\xi|^2}{2}.
\]
{Consequently (\ref{multiplier}) is satisfied}.

{Finally, if $0\neq \lambda = \pm \lambda_3$, we note that
\bFormula{d7a} \lambda_3 (|\xi|^2, k, \omega) = \omega |k|
\frac{1}{\lambda_1 (|\xi|^2, k, \omega)},\quad \lambda_3^2 - k^2 =
\frac{ k^2 (\omega^2 - \lambda_1^2) }{\lambda_1^2}, \eF
\[
\lambda^2_1 -\omega^2= \frac{ -\omega^2 + |\xi|^2 +|k|^2 + \sqrt{
(\omega^2 + |\xi|^2 + k^2 )^2 - 4 k^2 \omega^2} }{2}\ge |\xi|^2/2.
\]

Using identity (\ref{d7a})$_1$, we write $\vc E$ in terms of
$\lambda_1$ and verify (\ref{multiplier}) employing this explicit
formula and estimates (\ref{ev}), (\ref{d7a}).}

\subsection{Frequency cut-off}

As we shall see below, it will be convenient to approximate the
initial data for the problem (\ref{d3}), (\ref{d4}) by a frequency
truncation represented by a function $\psi \in \DC(0, \infty)$.

 Accordingly, solutions of the problem (\ref{d3}), (\ref{d4})
will be, in a certain way, composed of the quantities
\bFormula{d7} Z(t, x_h, k, \omega) = \mathcal{F}^{-1}_{\xi \to
x_h} \left[ \exp \Big( \pm {\rm i} \lambda_j ({|\xi|^2}, k, \omega)
t \Big) \psi (|\xi|) \hat h (\xi) \right] ,\  j = 1,2,3,4, \ k \in
Z, \ \omega \in (0, 1), \eF where $\hat h$ stands for the Fourier
transform of the ``initial data''. Our goal will be to derive
suitable \emph{dispersive estimates} for the mapping $h \mapsto
Z$.

We start with the $L^1 - L^\infty$ decay estimates. To this end,
write \bFormula{d8} \| Z(t, \cdot , k, \omega) \|_{L^\infty(R^2)}
\leq \left\| \mathcal{F}^{-1}_{\xi \to x_h} \left[ \exp \Big( \pm
{\rm i} \lambda_j ({ |\xi|^2 }, k, \omega) t \Big) \psi
(|\xi|)\right] \right\|_{L^\infty(R^2)} \| h \|_{L^1(R^2)}. \eF
Consequently,
\bFormula{d9} \mathcal{F}^{-1}_{\xi \to x_h} \left[ \exp \Big( \pm
{\rm i} \lambda_j ({|\xi|^2}, k, \omega) t \Big) \psi
(|\xi|)\right] (x_h) \eF
{
$$ ={ {\pi \sqrt 2}}\int_0^{2\pi}\int_0^\infty \exp
\Big( \pm {\rm i} \lambda_j (|\xi|^2, k, \omega) t \Big) \psi
(|\xi|)\exp\Big({\rm i} |\xi||x_h| \sin \theta\Big)  |\xi|\ {\rm
d}|\xi|{\rm d}\theta
$$
}
$$
= { {\pi \sqrt 2}}\int_0^\infty \exp \Big( \pm {\rm i} {
\lambda_j} (r^2, k, \omega) t \Big) \psi(r) r J_0(r |x_h|) \ {\rm
d}r,
$$
 where the symbols $J_m$, $m=0,1, \dots$ denote the Bessel
functions, cf. Guo, Peng, and Wang \cite{GuPeWa}.

Finally, performing a simple change of variables, we get
\bFormula{d10} \mathcal{F}^{-1}_{\xi \to x_h} \left[ \exp \Big(
\pm {\rm i} \lambda_j ({|\xi|^2}, k, \omega) t \Big) \psi
(|\xi|)\right] (x_h) = { \frac {\pi \sqrt 2}{2}} \int_0^\infty
\exp \Big( \pm {\rm i}{ \lambda_j} (z, k, \omega) t \Big)
\psi(\sqrt{z}) J_0(\sqrt{z}|x_h|) \ {\rm d}z. \eF

\subsection{Decay estimates}

Supposing $\lambda_j \ne 0$ we derive the desired decay estimates. To this end, we use van Corput's lemma, see Stein
\cite[Chapter 8.1.2, Proposition 2 and Corollary]{STEIN1}:

\Cbox{Cgrey}{

\bLemma{vC}
Let $\Lambda = \Lambda (z)$ be a smooth function away from the origin,
\[
\partial_z \Lambda (z) \ \mbox{monotone},
\
|\partial_z \Lambda (z)| \geq \Lambda_0 > 0
\]
for all $z \in [a,b]$, $0< a < b < \infty$. Let $\Phi$ be a smooth function on $[a,b]$.

Then
\[
\left| \int_{a}^b \exp \left( {\rm i} \Lambda(z) t \right) \Phi(z) \ {\rm d}z \right| \leq c
\frac{1}{t \Lambda_0} \left[ | \Phi(b)| + \int_a^b |\partial_z \Phi (z) | \ {\rm d}z   \right],
\]
where $c$ is an absolute constant independent of $\Lambda$ and $\Phi$.
\eL

}

Going back to the oscillatory integral (\ref{d10}), we distinguish two cases.

\subsubsection{Case $|x_h| > t^\beta$}

Using the decay properties of $J_0$ and the fact that $\psi$ is
compactly supported away from zero, we get \bFormula{d12} \left|
\mathcal{F}^{-1}_{\xi \to x_h} \left[ \exp \Big( \pm {\rm i}
\lambda_j ({|\xi|^2}, k, \omega) t \Big) \psi (|\xi|)\right] (x_h)
\right| \leq c(\psi) t^{-\beta/2} \ \mbox{whenever}\ |x_h| >
t^\beta. \eF

\subsubsection{Case $|x_h| \leq t^\beta $}

The idea is to use van Corput's lemma.
Let $[a,b]$ be a closed interval, $a > 0$, containing the support of $\psi (\sqrt{z})$. In accordance with the hypotheses of
Lemma \ref{LvC}, we have to verify that
\begin{itemize}
\item
\[
\left| \partial_z {\lambda_j} (z, k, \omega) \right| \geq
\Lambda_0 (a,b,\omega, k) > 0 \ \mbox{for all}\ z \in [a,b], \
\omega \in (0, 1) ;
\]
\item
$\partial_z \Lambda_j$ is monotone on $[a,b]$.
\end{itemize}

If the two conditions are satisfied, we get \bFormula{d13} \left|
\mathcal{F}^{-1}_{\xi \to x_h} \left[ \exp \Big( \pm {\rm i}
\lambda_j ({|\xi|^2}, k, \omega) t \Big) \psi (|\xi|)\right] (x_h)
\right|  \leq c(\psi) \frac{1}{\Lambda_0 (a,b, k,\omega) t} \Big(
1 + t^{\beta/2} \Big) \ \mbox{for}\ |x_h| \leq t^\beta,
\eF
 where again
we have used the properties of the Bessel functions, namely,
\[
{J'}_0(z) = - J_1(z).
\]

{Now, our goal is to verify the hypotheses of Lemma
\ref{LvC}}. We have
\[
\partial_z {\lambda_1}(z,k,\omega) = \frac{1}{2} \frac{1}{ \lambda_1(z,k,\omega) } \beta_1 (z,k,\omega),\
\beta_1 (z,k,\omega) = \frac{1}{2} \left( 1 + \frac{ \omega^2 + z + k^2 }{\sqrt{ (\omega^2 + z + k^2)^2 - 4 \omega^2 k^2 }} \right),
\]
and, similarly,
\[
\partial_z {\lambda_3}(z,k,\omega) = \frac{1}{2} \frac{1}{ \lambda_3 (z,k,\omega) } \beta_3 (z,k,\omega),\
\beta_3 (z,k,\omega) = \frac{1}{2} \left( 1 - \frac{ \omega^2 + z + k^2 }{\sqrt{ (\omega^2 + z + k^2)^2 - 4 \omega^2 k^2 }} \right), \ k \ne 0.
\]
Furthermore,
\[
\partial_z \beta_1 (z,k,\omega) = - \frac{2 \omega^2 k^2} {\left[  (\omega^2 + z + k^2)^2 - 4 \omega^2 k^2 \right]^{3/2} },
\
\partial_z \beta_3 (z,k,\omega) = \frac{2 \omega^2 k^2} {\left[  (\omega^2 + z + k^2)^2 - 4 \omega^2 k^2 \right]^{3/2} }.
\]

Summing up the previous relations, we conclude that
\[
\partial_z {\lambda_1} (z,k,\omega) \ \mbox{is a decreasing functions of}\ z ,
\]
while
\[
\partial_z {\lambda_3} (z,k,\omega) \ \mbox{is an increasing functions of}\ z \ \mbox{for}\ k \ne 0.
\]
Consequently, we deduce that
\bFormula{d14}
\partial_z {\lambda_1} (z,k,\omega) \geq \partial_z {\Lambda_1} (b,k,\omega) \geq \Lambda (\psi,k)  \ \mbox{for}\ z \in [a,b],\
\omega \in (0,1).
\eF

Finally, we have
\[
\partial_{z} {\lambda_3} (z,k,\omega) < 0, \ \left| \partial_{z} {\lambda_3} (z,k,\omega) \right|
\geq | \partial_z {\lambda_3} (b,k,\omega) |,
\]
where
\[
| \partial_z {\lambda_3} (b,k,\omega) | \geq \frac{1}{2}
\left( \frac{ \omega^2 + b + k^2}{ \sqrt{ (\omega^2 + b + k^2)^2 -
4 \omega^2 k^2 } } - 1 \right) \left( \frac{ \omega^2 + b + k^2 -
\sqrt{ (\omega^2 + b + k^2)^2 - 4 \omega^2 k^2 }  }{2}
\right)^{-1/2}
\]
\[
= \sqrt{2} \omega |k| \left( (\omega^2 + b + k^2 )^2 - 4 \omega^2 k^2 \right)^{-1/2} \left(
\omega^2 + b + k^2 + \sqrt{(\omega^2 + b + k^2)^2 - 4 \omega^2 k^2 } \right)^{-1/2}
\]
\[
\geq c(\psi,k) \omega \ \mbox{for}\ k \ne 0, \
\omega \in (0,1).
\]

Thus, reviewing (\ref{d12}), (\ref{d13}) we may infer that
\bFormula{d15} \left\| \mathcal{F}^{-1}_{\xi \to x_h} \left[ \exp
\Big( \pm {\rm i} \lambda_j ({|\xi|^2}, k, \omega) t \Big) \psi
(|\xi|)\right]  \right\|_{L^\infty(R^2)} \leq c (\psi, k) \max
\left\{  \frac{1}{\omega t^{1 - \beta/2}}  ; \frac{1}{t^{\beta/2}}
\right\},\ t>0, \eF
as soon as $\lambda_j \ne 0$, which gives rise to the decay
estimates \bFormula{d17} \| Z(t, \cdot , k, \omega)
\|_{L^\infty(R^2)} \leq c (\psi, k) \max \left\{ \frac{1}{\omega
t^{1 - \beta/2}}  ; \frac{1}{t^{\beta/2}} \right\} \| h
\|_{L^1(R^2)}. \eF
{Next, seeing that the mapping
$
h\mapsto \exp\Big({\rm i}{\lambda_j}(\xi,k,\omega) t\Big) h
$
is an isometry on $L^2(R^2)$, we deduce

\begin{equation}\label{d17+}
\| Z(t, \cdot , k, \omega) \|_{L^2(R^2)} \leq c\|h\|_{L^2(R^2)}
\end{equation}
}

Finally, interpolating (\ref{d17}) and (\ref{d17+}), we obtain the
$L^p - L^q$ estimates \bFormula{d18-} \| Z(t, \cdot , k, \omega)
\|_{L^p(R^2)} \leq   {c (\psi,p, k)} \max \left\{ \frac{1}{\omega
t^{1 - \beta/2}}  ; \frac{1}{t^{\beta/2}} \right\}^{1 -
\frac{2}{p}} \| h \|_{L^{p'}(R^2)} \ \mbox{for}\ p \geq 2, \
\frac{1}{p} + \frac{1}{p'} = 1,\ \lambda_j \ne 0. \eF

Keeping in mind that $\omega$ scales like $\ep^{m-1}$ while the
time $t$ is proportional to $\ep^{-m}$ we observe that taking
\[
0 < \beta < \frac{2}{m}
\]
yields the effective decay of $Z_\ep= Z(t/\ep^m,k,\omega)$ on any
compact subinterval of $(0,T]$. In particular, the optimal choice
$\beta=1/m$ gives rise to
\begin{equation}\label{d18}
\left\| Z\left(\frac{t}{\ep^m}, \cdot , k, \omega \right) \right\|_{L^p(R^2)}\le c \ \ep^{\frac{1}{2} - \frac{1}{p}}
\max \left\{ \frac{1}{t^{1 - 1/2m}}; \frac{1}{t^{1/2m}} \right\}^{1 - \frac{2}{p}} \|h\|_{L^{p'}(R^2)}, \,\;p\geq 2,\;\lambda_j\neq
0,\; t\in (0,T].
\end{equation}

\section{Convergence}
\label{c}

In this final part, we complete the proof of Theorem \ref{Tm1}. The basic idea is to use the relative entropy inequality (\ref{p6}) for a suitable choice of  test functions $r$ and $\vc{U}$.

\subsection{Initial data decomposition}

We start be introducing suitable smoothing operators imposed on the initial data. Taking a family of smooth functions
\[
\psi_\delta \in \DC(0, \infty), 0 \leq \psi_\delta \leq 1, \ \psi_\delta \nearrow 1 \ \mbox{as}\ \delta \to 0,
\]
and
\[
\phi_\delta = \phi_\delta(x_h) \in \DC(R^2),\ 0 \leq \phi_\delta \leq 1, \ \phi_\delta \nearrow 1 \ \mbox{as}\ \delta \to 0,
\]
we introduce
\bFormula{d20} \left[ \vr^{(1)}_0 \right]_{\delta}(x_h,x_3) =
\frac 1{\sqrt 2}\sum_{|k|\le 1/\delta} \mathcal{F}^{-1}_{\xi \to
x_h} \left[ \psi_\delta (|\xi|) \widehat { \left( \vr_{0}^{(1)}
\phi_\delta \right) }(\xi,k) \right] \exp \left( -{\rm i} {k x_3}
\right),  \eF
 and, similarly, \bFormula{d21} \left[ u_{0,j}
\right]_{\delta} (x_h,x_3) =\frac 1{\sqrt 2}\sum_{|k|\le 1/\delta}
\mathcal{F}^{-1}_{\xi \to x_h} \left[ \psi_\delta (|\xi|) \widehat {
\left( u_{0,j} \phi_\delta   \right) }(\xi,k) \right] \exp \left(-
{\rm i} {k x_3} \right), \ j=1,2,3. \eF

Now, we write the initial data in the form
\bFormula{d22}
\left[ \vr^{(1)}_0 \right]_{\delta} = s_{0, \ep, \delta} + q_{0,\ep, \delta} ,\ \mbox{where} - \Delta_h q_{0, \ep, \delta} + \omega^2 q_{0, \delta} = \omega \int_0^1 {\rm curl}_h \left[ [\vc{u}_0]_h \right]_\delta \ {\rm d}x_3 + \omega^2
\int_0^1 \left[ \vr^{(1)}_0 \right]_{\delta} {\rm d}x_3
\eF
\bFormula{d23}
\left[ \vc{u}_0 \right]_{\delta} = \vc{V}_{0, \ep, \delta} + \vc{v}_{0, \ep, \delta} , \ \mbox{with}\ \omega [v_{0,\ep, \delta}]_1 = - \partial_{x_2} q_{0, \ep, \delta},\ \omega [v_{0,\ep, \delta}]_2 =  \partial_{x_1} q_{0, \ep, \delta}.
\eF

Finally, we choose the functions $r$, $\vc{U}$ in the relative entropy inequality as
\bFormula{d23a}
r = r_{\ep,\delta} = \tilde \vre + \ep^m (q_{\ep, \delta} + s_{\ep, \delta} ),\ \vc{U} = \vc{U}_{\ep, \delta} = \vc{v}_{\ep, \delta} + \vc{V}_{\ep, \delta},
\eF
where $[s_{\ep, \delta}, \vc{V}_{\ep, \delta}]$ is the unique solution of the acoustic-Rossby system (\ref{d1}), (\ref{d2}), emanating from the initial data
\[
s_{\ep, \delta}(0, \cdot) = s_{0, \ep, \delta},\ \vc{V}_{\ep, \delta}(0, \cdot) = \vc{V}_{0,\ep, \delta},
\]
while the functions $q_{\ep, \delta}$, $\vc{v}_{\ep, \delta}$ are interrelated through
\begin{equation}\label{d23b}
\omega \vc{f} \times \vc{v}_{\ep, \delta} + \Grad q_{\ep, \delta} = 0,
\end{equation}
where $q_{\ep, \delta}$ is the unique solution of the problem
\bFormula{d24}
\partial_t \left( \Delta_h q_{\ep, \delta} - \omega^2 q_{\ep, \delta} \right) + \frac{1}{\omega} {\nabla_h^{\perp}} q_{\ep,\delta} \cdot \Grad \left( \Delta_h q_{\ep, \delta} - \omega^2 q_{\ep, \delta} \right) = 0, \ q_{\ep, \delta}(0, \cdot) = q_{0, \ep, \delta}.
\eF

\subsection{Decay of the oscillatory component}

First we claim that, in view of the dispersive estimates (\ref{d18}) (with $0 < \beta < 2/m$), and the properties of the eigenvectors
of the matrix $\mathcal{A}$, discussed in detail in Section \ref{EV},
we get
\bFormula{d25} s_{\ep, \delta} \to 0 \ \mbox{in}\ L^p(0,T; W^{l,
q}(\Omega)), \vc{V}_{\ep, \delta} \to 0 \ \mbox{in}\ L^p(0,T;
W^{l, q}(\Omega)) \ \mbox{as} \ \ep \to 0 \eF
for any \emph{fixed} $\delta > 0$, $1 \leq p < \infty$, $2 < q \leq \infty$ and
$l=0,1,\dots$

We note that
\[
\left[ \begin{array}{c} s_{\ep, \delta}(t,x_h,k) \\ \\ \vc{V}_{\ep, \delta}(t, x_h, k \end{array} \right] =
\mathcal{F}^{-1}_{\xi \to x_h} \left[ \mathcal{Q}^T (\xi, \omega ,k)
\left[ \begin{array}{c} \exp \left( {\rm i} \lambda_1 (\xi, \omega,k) \frac{t}{\ep^m} \right) , 0 , 0 , 0
\\ \\ 0, \exp \left( {\rm i} \lambda_2 (\xi, \omega,k) \frac{t}{\ep^m} \right), 0 , 0 \\ \\
0, 0, \exp \left( {\rm i} \lambda_3 (\xi, \omega,k) \frac{t}{\ep^m} \right), 0 \\ \\
0, 0, 0, \exp \left( {\rm i} \lambda_4 (\xi, \omega,k) \frac{t}{\ep^m} \right) \end{array} \right] \mathcal{Q}(\xi, \omega, k) \psi(|\xi|) \hat h_0(\xi) \right]
\]
\[
=\mathcal{F}^{-1}_{\xi \to x_h} \left[ \tilde \psi (|\xi|) \mathcal{Q}^T (\xi, \omega ,k)
\left[ \begin{array}{c} \exp \left( {\rm i} \lambda_1 (\xi, \omega,k) \frac{t}{\ep^m} \right) , 0 , 0 , 0
\\ \\ 0, \exp \left( {\rm i} \lambda_2 (\xi, \omega,k) \frac{t}{\ep^m} \right), 0 , 0 \\ \\
0, 0, \exp \left( {\rm i} \lambda_3 (\xi, \omega,k) \frac{t}{\ep^m} \right), 0 \\ \\
0, 0, 0, \exp \left( {\rm i} \lambda_4 (\xi, \omega,k) \frac{t}{\ep^m} \right) \end{array} \right] \tilde \psi(|\xi|)
\mathcal{Q}(\xi, \omega, k) \psi(|\xi|) \hat h_0(\xi) \right], \ \omega = \ep^{m-1},
\]
where $\tilde \psi \in \DC(0, \infty)$ has been chosen so that $\psi \tilde \psi = \psi$.

In accordance with Section \ref{EV}, the quantities
$
\tilde \psi \mathcal{Q}, \tilde \psi \mathcal{Q}^T
$
are $L^p-$ Fourier multipliers with norm independent of $\omega$. Thus the desired decay estimates (\ref{d25}) follow from (\ref{d18}).

\subsection{Convergence of the non-oscillatory component}

We introduce the scaled function $\tilde q _{\ep, \delta} = q_{\ep, \delta}/\omega$ and observe that
\bFormula{d26}
\partial_t \left( \Delta_h \tilde q_{\ep, \delta} - \omega^2 \tilde q_{\ep, \delta} \right)
+  { \nabla_h^{\perp}} \tilde q_{\ep,\delta} \cdot \nabla_h
\left( \Delta_h \tilde q_{\ep, \delta} - \omega^2 \tilde q_{\ep,
\delta} \right), \ { q_{\ep, \delta}(0) = q_{0,\ep, \delta}},
\eF
\[
\Delta_h \tilde q_{0, \ep, \delta} - \omega^2 \tilde q_{0, \delta} = \int_0^1 {\rm curl}_h \left[ [\vc{u}_0]_h \right]_\delta \ {\rm d}x_3 + \omega
\int_0^1 \left[ \vr^{(1)}_0 \right]_{\delta} {\rm d}x_3.
\]

Since the initial data are regular, we may use the result of Oliver \cite[Theorem 3]{Oli} to deduce that
\bFormula{d27}
\{ \Delta \tilde q_{\ep, \delta} + \ep^2 \tilde q_{\ep, \delta} \}_{\ep > 0}
\ \mbox{is bounded in}\ C^r([0,T]; W^{l,2}(\Omega)),\ r\geq 0, \ l=0,1,\dots
\eF
Moreover, as $\Delta \tilde q_{\ep, \delta}$ satisfies that transport equation (\ref{d26}) with the initial data in $L^p(R^2)$,
we get
\bFormula{d28}
\{ \Delta \tilde q_{\ep, \delta} + \omega^2 \tilde q_{\ep, \delta} \}_{\ep > 0}
\ \mbox{is bounded in}\ L^\infty([0,T]; L^p(\Omega)) \ \mbox{for any}\ 1 < p < \infty.
\eF

Next, we recall the ``energy estimates'' that can be obtained multiplying (\ref{d26}) on $\tilde q_{\ep, \delta}$ and integrating by parts:
\bFormula{d29}
\int_{R^2} \left(  |\Grad \tilde q_{\ep, \delta} |^2 + \omega^2 | \tilde q_{\ep, \delta}|^2 \right)(\tau, \cdot) \ {\rm d}x =
\int_{R^2} \left(  |\Grad \tilde q_{0,\ep, \delta} |^2 + \omega^2 |\tilde q_{0,\ep, \delta}|^2 \right) \ {\rm d}x.
\eF
Note that
\[
\int_{R^2}  {\nabla_h^{\perp}} \tilde q_{\ep, \delta} \cdot
\nabla_h \Delta_h \tilde q_{\ep, \delta} \tilde q_{\ep, \delta}
\dx = - \int_{R^2}  {\nabla_h^{\perp}} \tilde q_{\ep, \delta}
\cdot \nabla_h \tilde q_{\ep, \delta} \Delta_h \tilde q_{\ep,
\delta} \dx = 0.
\]

Since
\[
\vc{v}_{\ep, \delta} = {\nabla_h^{\perp}} \tilde q_{\ep,
\delta},
\]
we get
\bFormula{d30}
\{ \vc{v}_{\ep, \delta} \}_{\ep > 0}
\ \mbox{is bounded in}\ C^r([0,T]; W^{l,2}(R^2)),\ r\geq 0, \ l=0,1,\dots
\eF

Finally, we compute
\[
\partial_t \tilde q_{\ep, \delta} = (\Delta_h - \omega^2)^{-1} \left[ \nabla_h \left( \vc{v}_{\ep, \delta} {\rm curl} \vc{v}_{\ep, \delta} \right) \right],
\]
where, furthermore,
\[
\vc{v} {\rm curl}_h \vc{v} = \Big[ v_1 \left( \partial_{x_1} v_2 - \partial_{x_2} v_1 \right) ; v_2 \left( \partial_{x_1} v_2 - \partial_{x_2} v_1 \right) \Big],
\]
with
\[
v_1 \partial_{x_1} v_2  = \partial_{x_1}(v_1 v_2) + \frac{1}{2} \partial_{x_2} v_2^2,\
v_2 \partial_{x_2} v_1 = \partial_{x_2}(v_1 v_2) + \frac{1}{2} \partial_{x_1} v_1^2.
\]

Consequently, we may infer that
\bFormula{d30a}
\{ \partial_t \tilde q_{\ep, \delta} \}_{\ep > 0} \ \mbox{is bounded in}\ C^r([0,T]; W^{l,q}(R^2)),\ r\geq 0,\ q > 1, \ l=0,1,\dots
\eF
All the above estimates may depend on $\delta$ but are uniform with respect to $\ep \searrow 0$.

In view of the above estimates, it is easy to pass to the limit
for $\ep \to 0$ in order to get
\bFormula{A2}\vc{v}_{\ep, \delta} \to \vc{v}_{\delta}, \
\partial_t \vc{v}_{\ep, \delta} \to \partial_t \vc{v}_\delta \
\mbox{weakly-(*) in}\ L^\infty(0,T; W^{l,2}(R^2)), \ l=0,1,\dots,
\eF where
\bFormula{d31} \vc{v}_\delta , \
\partial_t \vc{v}_\delta \in C([0,T]; W^{l,2}(R^2)), \
l=0,1,2,\dots. \eF

{We have, in particular,
$$
 \vc{v}_{\ep, \delta} \to \vc{v}_{\delta}\;\mbox{in $L^q(0,T;
L^q_{\rm loc}(\Omega))$, $1\le q<\infty$};
$$
whence,  by virtue of (\ref{d26}),} {
$$
\partial_t {\rm curl}_h\vc v_\delta+\vc v_\delta\cdot\nabla_h{\rm
curl}_h\vc v_\delta=0.
$$
Seeing that
$$
{\rm curl}_h{\rm div}_h(\vc v_\delta\otimes\vc v_\delta)={\rm
curl}_h(\vc v_\delta\cdot\nabla_h\vc v_\delta)= \vc
v_\delta\cdot\nabla_h{\rm curl}_h\vc v_\delta,
$$
we deduce existence of \bFormula{d31+} \Pi_\delta\in C([0,T];
W^{l,2}(R^2)), \ l=0,1,2,\dots, \eF
where the couple $(\vc
v_\delta, \Pi_\delta)$ is
the unique solution
of the Euler system (\ref{i9}), emanating from the initial data
\bFormula{d32} \vc{v}_{\delta}(0, \cdot) = \int_0^1 \vc{H}_h
\left[ \left[ \vc{u}_0 \right]_{\delta} \right] \ {\rm d}x_3. \eF
}

\subsection{Relative entropy inequality}

We return to the relative entropy inequality (\ref{p6}), where $\vr = \vre$, $\vu = \vue$ and the test functions $r$ and $\vc{U}$ are given by
(\ref{d23a}). In what follows, we examine step by step all terms on the right-hand side of (\ref{p6}) and perform the limits; first for $\ep \to 0$, then for $\delta \to 0$.

\subsubsection{Initial data}

We have
\bFormula{r5a}
\mathcal{E}_\ep \left( \vr_{0,\ep}, \vu_{0,\ep} \ \Big| \ r(0,\cdot) , \vc{U}(0,\cdot) \right)
\eF
\[
=
\intO{ \frac{1}{2} \vr_{0,\ep} |\vu_{0,\ep} - [\vu_0]_\delta  |^2 }
\]
\[
+  \intO{ \left[ \frac{1}{\ep^{2m}}
\left( H \left( 1 + \ep^m \vr^{(1)}_{0,\ep} \right) -  H'(1 + \ep^m [ \vr^{(1)}_0 ]_{\delta} ) \left((\vr^{(1)}_{0,\ep} - [\vr^{(1)}_0] _\delta
\right) - H (1 + \ep^m [ \vr^{(1)}_0 ]_{\delta} )\right) \right] }
\]
\[
\leq c \left( \| \vu_{0,\ep} - [ \vu_{0} ]_{\delta} \|^2_{L^2(\Omega;R^3)} +
\| \vr^{(1)}_{0,\ep} - [ \vr^{(1)}_{0} ]_{\delta} \|^2_{L^2(\Omega;R^3)} \right)
\]
\[
\to c \left( \| \vu_{0} - [ \vu_{0} ]_{\delta} \|^2_{L^2(\Omega;R^3)} +
\| \vr^{(1)}_{0} - [ \vr^{(1)}_{0} ]_{\delta} \|^2_{L^2(\Omega;R^3)} \right) \ \mbox{as} \ \ep \to 0.
\]
The most left quantity obviously tends to zero as $\delta \searrow 0$.

\subsubsection{Viscosity}

We write \bFormula{r5b+} \ep^\alpha \left| \int_0^\tau \intO{
\tn{S} (\Grad \vc{U}_{\ep, \delta} )  : \Grad (\vc{U}_{\ep,\delta}
- \vu_\ep) } \ \dt \right| \leq \ep^\alpha c_1(\delta) \int_0^\tau
\intO{ \left| \Grad (\vc{U}_{\ep,\delta} - \vu_\ep) \right| } \
\dt, \eF
and, by Korn's inequality,
\[
\ep^\alpha \int_0^\tau \intO{ \left| \Grad (\vc{U}_{\ep,\delta} - \vu_\ep) \right| } \ \dt \leq
\frac{\ep^\alpha}{2} \int_0^\tau \intO{ \left( \tn{S}_\ep(\Grad \vue) - \tn{S}_\ep (\Grad \vc{U}_{\ep, \delta} \right):
\Grad \left( \vc{U}_{\ep, \delta} - \vue \right) } \ \dt  + c_2 \ep^\alpha.
\]

\subsubsection{Forcing term}\label{6.4.3}

We have
\bFormula{r11} \frac 1{\ep^{2m}}\Big[ (r_{\ep,\delta} - \vre)
\partial_t H'(r_{\ep, \delta}) +{\Grad  \Big(H'(r_{\ep,
\delta})-H'(\tilde\vre)\Big)} \cdot (r_{\ep, \delta} \vc{U}_{\ep,
\delta} - \vre \vue ) \Big] \eF
\[
+ \frac 1{\ep^{2m}}\Div \vc{U}_{\ep, \delta} \Big( p(r_{\ep,
\delta})-p(\vre) \Big)+ \frac
1{\ep^{2n}}(r_{\ep,\delta}-\vre)\Grad G\cdot\vc U_{\ep,\delta}\]
\[
= \frac 1{\ep^{2m}}\Big[ p(r_{\ep, \delta})  - p'(r_{\ep, \delta})
(r_{\ep, \delta} - \vre) - p(\vre) \Big] \Div \vc{U}_{\ep, \delta}
+ \frac 1{\ep^{2m}}(r_{\ep, \delta} - \vre) H''(r_{\ep, \delta})
\Big[
\partial_t r_{\ep, \delta} + \Div (r_{\ep, \delta} \vc{U}_{\ep,
\delta}) \Big]
\]
\[
+\frac 1{\ep^{2m}}\Big[\vre \Grad H'(r_{\ep, \delta}) \cdot (
\vc{U}_{\ep, \delta} - \vue) {-\Grad H'(\tilde\vre)\cdot (
r_{\ep,\delta}\vc{U}_{\ep, \delta} - \vre\vue)}\Big]+ \frac
1{\ep^{2n}}(r_{\ep,\delta}-\vre)\Grad G\cdot\vc U_{\ep,\delta},
\]
where \bFormula{rr11-}
\partial_t r_{\ep, \delta} + \Div (r_{\ep, \delta} \vc{U}_{\ep, \delta}) = (\tilde \vre - 1) \Div \vc{V}_{\ep,\delta}
+ \Grad \tilde \vre \cdot \vc{U}_{\ep, \delta} + \eF
\[ \ep^m
\partial_t q_{\ep, \delta} + \ep^m \Div \left( (s_{\ep, \delta} +
q_{\ep, \delta})(\vc{v}_{\ep, \delta} + \vc{V}_{\ep, \delta} )
\right).
\]
and {
\begin{equation}\label{rr11+}
\vre \Grad H'(r_{\ep, \delta}) \cdot ( \vc{U}_{\ep, \delta} -
\vue)-\Grad H'(\tilde\vre)\cdot(r_{\ep,\delta}\vc
U_{\ep,\delta}-\vre\vue)
\end{equation}
\[
=\vre \Grad \Big[ H'(r_{\ep, \delta}) - H''(\tilde \vre)  (r_{\ep,
\delta} - \tilde \vre) - H'(\tilde \vre) \Big]\cdot ( \vc{U}_{\ep,
\delta} - \vue)
\]
\[
+ \vre \Grad H''(\tilde \vre) \cdot (\vc{U}_{\ep, \delta} - \vue)
(r_{\ep, \delta} - \tilde \vre) + \vre \ep^m H''(\tilde \vre)
\Big( \Grad s_{\ep, \delta} + \Grad q_{\ep, \delta} \Big) \cdot (
\vc{U}_{\ep, \delta} - \vue)
\]
\[
+(\vre- r_{\ep,\delta})\Grad H'(\tilde\vre)\cdot\vc U_{\ep,\delta}
.
\]
Here we have used the fact that the term  $\vre \Grad H'(\tilde
\vre) \cdot ( \vc{U}_{\ep, \delta}-\vue)$ in the
expansion of $\vre \Grad H'(r_{\ep, \delta}) \cdot ( \vc{U}_{\ep,
\delta}-\vue)$ cancels with the same term in the
expansion of $\Grad H'(\tilde\vre)\cdot(r_{\ep,\delta}\vc
U_{\ep,\delta}-\vre\vue)$. }
%
%

{ Now we use formula (\ref{r11}) with the second and third
terms at the right hand side expressed through formulas
(\ref{rr11-}) and (\ref{rr11+}): we keep the terms $\frac
1{\ep^m}\vre H''(\tilde \vre) \Big( \Grad s_{\ep, \delta} + \Grad
q_{\ep, \delta} \Big) \cdot ( \vc{U}_{\ep, \delta} - \vue)$ and $
\frac 1{\ep^{m}}(r_{\ep,\delta}-\vre)\Grad
H'(r_{\ep,\delta})\left[\partial_t q_{\ep, \delta} + \Div \left(
(s_{\ep, \delta} + q_{\ep, \delta})(\vc{v}_{\ep, \delta} +
\vc{V}_{\ep, \delta} ) \right)\right]$  as they are, and estimate
the decay of all remaining terms as $\ep\to 0$.}

{We observe that, thanks to (\ref{d26}), (\ref{d29}),
\bFormula{R12} \sup_{t \in (0,T)} \| q_{\ep, \delta}(t, \cdot)
\|_{L^2(R^2)} \leq c,\; \sup_{t \in (0,T)} \| \Delta q_{\ep,
\delta} (t, \cdot) \|_{L^2(R^2)} \leq c\ep^{m-1},\; \sup_{t \in
(0,T)} \| \partial_t q_{\ep, \delta} (t, \cdot) \|_{L^2(R^2)} \leq
c\ep^{m-1}. \eF Now, we use (\ref{p2}--\ref{p3}),
(\ref{rr2}--\ref{rr5}), (\ref{d23a}), (\ref{d25}), (\ref{d30}),
(\ref{R12}) to deduce the following estimates:
{\bFormula{1} \frac
1{\ep^{2n}}\Big|\int_0^\tau\intO{(r_{\ep,\delta}-\vre)\Grad
G\cdot\vc U_{\ep,\delta}}{\rm d}t\Big|\le c\ep^{m-2n}, \eF }
%
\medskip
\bFormula{2} \frac 1{\ep^{2m}}\int_0^\tau\intO{\Big[ p(r_{\ep,
\delta}) - p'(r_{\ep, \delta}) (r_{\ep, \delta} - \vre) - p(\vre)
\Big] \Div \vc{U}_{\ep, \delta}}{\rm d}t \eF
$$
= \frac 1{\ep^{2m}}\int_0^\tau\intO{\Big[ p(r_{\ep, \delta})  -
p'(r_{\ep, \delta}) (r_{\ep, \delta} - \vre) - p(\vre) \Big]_{\rm
res} \Div \vc{V}_{\ep, \delta}}{\rm d}t
$$
$$
+\frac 1{\ep^{2m}}\int_0^\tau\intO{\Big[ p(r_{\ep, \delta})  -
p'(r_{\ep, \delta}) (r_{\ep, \delta} - \vre) - p(\vre) \Big]_{\rm
ess} \Div \vc{V}_{\ep, \delta}}{\rm d}t= h(\ep,\delta),
$$
\medskip
\bFormula{3} \frac 1{\ep^{2m}}\Big|\int_0^\tau\intO{(r_{\ep,
\delta} - \vre) H''(r_{\ep, \delta})\Big((\tilde \vre - 1) \Div
\vc{V}_{\ep,\delta} + \Grad \tilde \vre \cdot \vc{U}_{\ep,
\delta}\Big)}{\rm d} t\Big| \le c\ep^{m-2n}, \eF
\medskip
\bFormula{4} \frac 1{\ep^{2m}}\int_0^\tau\intO{\vre \Grad \Big[
H'(r_{\ep, \delta}) - H''(\tilde \vre) (r_{\ep, \delta} - \tilde
\vre) - H'(\tilde \vre) \Big]\cdot ( \vc{U}_{\ep, \delta} -
\vue)}{\rm d} t= h(\ep,\delta), \eF
\medskip
\bFormula{5} \frac 1{\ep^{2m}}\Big(\Big|\intO{\vre \Grad
H''(\tilde \vre) \cdot (\vc{U}_{\ep, \delta} - \vue) (r_{\ep,
\delta} - \tilde \vre)}{\rm d}t\Big| +\Big|\int_0^\tau\intO{(\vre-
r_{\ep,\delta})\Grad H'(\tilde\vre)\cdot\vc U_{\ep,\delta}}{\rm
d}t\Big|
\le c\ep ^{m-2n}. \eF

Here and hereafter, $h(\ep, \delta)$ denotes a generic function
having the property \bFormula{h} h(\ep, \delta) \to \tilde
h(\delta) \ \mbox{as}\ \ep \to 0, \ \tilde h(\delta) \to 0 \
\mbox{as} \ \delta \to 0. \eF } Note that the
dispersive decay estimates (\ref{d25}) play a crucial role in the analysis.

{Taking into account (\ref{r5a}), (\ref{r5b+}),  using the
identity (\ref{r11}) with the third and fourth terms expressed
through (\ref{rr11-}--\ref{rr11+}), and employing the asymptotic
behavior from formulas (\ref{1}--\ref{5}), we may rewrite the
relative entropy inequality (\ref{p6}) in the form \bFormula{R13}
\mathcal{E}_\ep \left( \vre, \vue \ \Big| \ r_{\ep, \delta},
\vc{U}_{\ep, \delta} \right) (\tau) \leq \int_0^\tau \intO{  \vre
\left( \partial_t \vc{U}_{\ep,\delta} + \vue \cdot \Grad
\vc{U}_{\ep, \delta} \right) \cdot \left( \vc{U}_{\ep, \delta} -
\vue \right) } \ \dt \eF
\[
 + \frac{1}{\ep} \int_0^\tau \intO{
\vre (\vc{f} \times \vue ) \cdot (\vc{U}_{\ep, \delta} - \vue) } \
\dt\]
\[
+ \frac{1}{\ep^{m}} \int_0^\tau \intO{ (r_{\ep,\delta}-\vre)
H''(r_{\ep,\delta}) \Big[
\partial_t q_{\ep, \delta} +  \Div \left( (s_{\ep, \delta} +
q_{\ep, \delta})\vc{U}_{\ep, \delta}  \right) \Big]    } \ \dt
\]
\[
+ \frac{1}{\ep^{m}} \int_0^\tau \intO{  \vre  H''(\tilde \vre)
\Big( \Grad s_{\ep, \delta} + \Grad q_{\ep, \delta} \Big) \cdot (
\vc{U}_{\ep, \delta} - \vue)  } \ \dt + h(\ep, \delta).
\]

\subsubsection{Coriolis force}
{We may write
\begin{equation}\label{A1}
\vre  H''(\tilde \vre) \Big( \Grad s_{\ep, \delta} + \Grad q_{\ep,
\delta} \Big) \cdot ( \vc{U}_{\ep, \delta} - \vue)
\end{equation}
\[
= \vre  (H''(\tilde \vre) - H''(1)) \Big( \Grad s_{\ep, \delta} +
\Grad q_{\ep, \delta} \Big) \cdot ( \vc{U}_{\ep, \delta} - \vue) +
\vre  H''(1) \Big( \Grad s_{\ep, \delta} + \Grad q_{\ep, \delta}
\Big) \cdot ( \vc{U}_{\ep, \delta} - \vue),
\]
where, by the same reasoning as in estimates (\ref{1}--\ref{5}),
$$
\frac 1{\ep^{2m}}\Big|\int_0^\tau\intO{\vre  (H''(\tilde \vre) -
H''(1)) \Big( \Grad s_{\ep, \delta} + \Grad q_{\ep, \delta} \Big)
\cdot ( \vc{U}_{\ep, \delta} - \vue)}{\rm d} t\Big|\le c\ep^{m-2n}
$$

Recalling our convention $H''(1) = 1$ we get for the second term
in (\ref{A1}),
\[
\frac 1{\ep^{m}}\vre  H''(1) \Big( \Grad s_{\ep, \delta} + \Grad
q_{\ep, \delta} \Big) \cdot ( \vc{U}_{\ep, \delta} - \vue) =
\]
\[
\vre  \Big( -  \partial_t \vc{V}_{\ep, \delta} - \frac 1\ep \vc{f}
\times \vc{V}_{\ep, \delta} - \frac 1\ep \vc{f} \times
\vc{v}_{\ep, \delta} \Big)\cdot ( \vc{U}_{\ep, \delta} - \vue)
\]
$$
=-\partial_t \vc V_{\ep,\delta}-\frac 1 {\ep}(\vc f\times\vc
U_{\ep,\delta})\cdot\vc u_\ep,
$$
where the last term cancels out with the Coriolis force term
$
\frac 1\ep\vc f\times\vu_\ep\cdot (\vc U_{\ep,\delta}-\vue)=\frac
1\ep(\vc f \times\vu_\ep)\cdot\vc U_{\ep,\delta}.
$

Thanks to formula (\ref{A1}) the relative entropy inequality
(\ref{R13}) reduces to \bFormula{R14} \mathcal{E}_\ep \left( \vre,
\vue \ \Big| \ r_{\ep, \delta}, \vc{U}_{\ep, \delta} \right)
(\tau) \leq \int_0^\tau \intO{  \vre \left( \partial_t
\vc{v}_{\ep,\delta} + \vue \cdot \Grad \vc{U}_{\ep, \delta}
\right) \cdot \left( \vc{U}_{\ep, \delta} - \vue \right) } \ \dt
\eF
\[
+ \frac{1}{\ep^{m}} \int_0^\tau \intO{ (r_{\ep, \delta} - \vre) H''(r_{\ep, \delta}) \Big[\partial_t q_{\ep, \delta} + \Div \left( (s_{\ep, \delta} + q_{\ep, \delta})\vc{U}_{\ep, \delta}  \right)                                 \Big]    } \ \dt
+ h(\ep, \delta).
\]

\subsubsection{Estimating the remaining terms}

We observe that
$$
\int_0^\tau \intO{  \vre \left( \partial_t \vc{v}_{\ep,\delta} +
\vue \cdot \Grad \vc{U}_{\ep, \delta} \right) \cdot \left(
\vc{U}_{\ep, \delta} - \vue \right) } \ \dt= \int_0^\tau \intO{
\vre\left(\vc v_{\ep,\delta}+  \vc U_{\ep,\delta} \cdot \Grad
\vc{U}_{\ep, \delta} \right) \cdot \left( \vc{U}_{\ep, \delta} -
\vue \right) } \ \dt
$$
$$
-\int_0^\tau \intO{  \vre (\vc U_{\ep,\delta}-\vue) \cdot \Grad
\vc{U}_{\ep, \delta}  \cdot \left( \vc{U}_{\ep, \delta} - \vue
\right) } \ \dt,
$$
where the terms in the first expression at the right hand side
containing the quantities $s_{\ep,\delta}$, $\vc V_{\ep,\delta}$
tend to $0$ in the limit $\lim_{\delta\to 0}\lim_{\ep\to 0}$
thanks to dispersive estimates (\ref{d25}). Consequently,
employing (\ref{p3}), we obtain
\[
\int_0^\tau \intO{  \vre \left( \partial_t \vc{v}_{\ep,\delta} +
\vue \cdot \Grad \vc{U}_{\ep, \delta} \right) \cdot \left(
\vc{U}_{\ep, \delta} - \vue \right) } \ \dt
\]
\[
\le \int_0^\tau \intO{   \left( \partial_t \vc{v}_{\ep,\delta} +
\vc{v}_{\ep, \delta} \cdot \Grad \vc{v}_{\ep, \delta} \right)
\cdot \left( \vc{v}_{\ep, \delta} - \vue \right) } \ \dt  + c
\int_0^\tau \mathcal{E}_\ep \left(\vre, \vue \ \Big| \ r_{\ep,
\delta}, \vc{U}_{\ep, \delta} \right)  \ \dt + h(\ep, \delta).
\]

Similarly, we deduce
\[
\frac{1}{\ep^{m}} \int_0^\tau \intO{ (r_{\ep, \delta} - \vre) H''(r_{\ep, \delta}) \Big[\partial_t q_{\ep, \delta} + \Div \left( (s_{\ep, \delta} + q_{\ep, \delta})\vc{U}_{\ep, \delta}  \right) \Big]    } \ \dt
\]
\[
=  \int_0^\tau \intO{
q_{\ep, \delta}
\partial_t q_{\ep, \delta} } \ \dt +
h(\ep, \delta),
\]
where we have used (\ref{rr2}--\ref{rr4}),(\ref{d23b}),
(\ref{R12}).

Now, we employ the energy equality (\ref{d29}) to observe that
\[
\int_0^\tau \intO{ \Big( \vc{v}_{\ep,\delta} \cdot \partial_t \vc{v}_{\ep, \delta} + q_{\delta, \ep} \cdot \partial_t q_{\ep, \delta} \Big) } \ \dt = 0.
\]
Consequently,
$$
\int_0^\tau \intO{   \left( \partial_t \vc{v}_{\ep,\delta} +
\vc{v}_{\ep, \delta} \cdot \Grad \vc{v}_{\ep, \delta}
\right)\cdot(\vc v_{\ep,\delta}-\vue)}{\rm d} t
$$
$$
+ \int_0^\tau \intO{ 
q_{\ep, \delta}
\partial_t q_{\ep, \delta} } \ \dt=-\int_0^\tau \intO{   \left( \partial_t \vc{v}_{\ep,\delta} +
\vc{v}_{\ep, \delta} \cdot \Grad \vc{v}_{\ep, \delta}
\right)\cdot\vue}{\rm d}t,
$$
where we have used the identity $\int_0^\tau \intO{ \vc{v}_{\ep,
\delta} \cdot \Grad \vc{v}_{\ep, \delta} \cdot\vc
v_{\ep,\delta}}{\rm d}t=0$.

We remark that
\[
\int_0^\tau \intO{   \left( \partial_t \vc{v}_{\ep,\delta} +
\vc{v}_{\ep, \delta} \cdot \Grad \vc{v}_{\ep, \delta} \right)
\cdot  \vue  } \ \dt \to \int_0^\tau \intO{   \left(
\partial_t \vc{v}_{\delta} + \vc{v}_{\delta} \cdot \Grad
\vc{v}_{\delta} \right) \cdot  \vu  } \ \dt
\]
\[
= \int_0^\tau \intO{ \Grad \Pi_\delta \cdot\vu } \ \dt = 0,
\]
where we have used (\ref{e3}--\ref{e4}), (\ref{A2}) and equations
(\ref{i9}) with $(\vc v,\Pi)$ replaced by $(\vc
v_\delta,\Pi_\delta)$. Consequently, the entropy inequality
(\ref{R14}) takes the form
$$
\mathcal{E}_\ep \left( \vre, \vue \ \Big| \ r_{\ep, \delta},
\vc{U}_{\ep, \delta} \right) (\tau) \leq h(\ep,\delta),
$$
where the function $h$ satisfies (\ref{h}). This finishes the
proof of Theorem \ref{Tm1}

\def\cprime{$'$} \def\ocirc#1{\ifmmode\setbox0=\hbox{$#1$}\dimen0=\ht0
  \advance\dimen0 by1pt\rlap{\hbox to\wd0{\hss\raise\dimen0
  \hbox{\hskip.2em$\scriptscriptstyle\circ$}\hss}}#1\else {\accent"17 #1}\fi}


\end{document}